\documentclass[11pt,reqno]{amsart} 
\usepackage{amsthm}
\usepackage{amssymb}
\usepackage{graphicx}
\usepackage{parskip}

\newtheorem{assert}{Assertion}
\newtheorem{prop1}{Property}

\newcommand{\comm}[1]{}
\renewcommand{\Re}{\operatorname{Re}}
\renewcommand{\Im}{\operatorname{Im}}
 
\newtheoremstyle{propstyle}
{15pt}
{1pt}
{\itshape}
{}
{\bfseries}
{.}
{.5em}
{}
\theoremstyle{propstyle}\newtheorem{prop}[prop1]{Property}

\graphicspath{%
    {converted_graphics/}
    {/}
}
\begin{document}

\title{A Reformulation of the Xi Function}


\author{Jon Breslaw}

\thanks{{\it Address:~} Department of Economics, Concordia University, Montreal, QC, Canada H3G 1M8}
\thanks{{\it Email:~}jon.breslaw@concordia.ca}

%



\date{December 4, 2015.}


\begin{abstract}
This paper proposes a reformulation of the Riemann Xi function in order to investigate its properties.  The reformulated function, which depicts the Xi function as the weighted sum of incomplete gamma functions, is validated, and a number of properties are established.  These properties are then used to analyze the intersection of the real and imaginary zero level curves. It is shown that a  pair of zeros off the critical line is not consistent with these
properties, thus validating the Riemann hypothesis.
\end{abstract}

\maketitle

\section{Introduction}

Riemann's original $\xi$ function is defined [3] for complex $z = x + i \, y$   as:
\begin{equation}
\xi(z)	=	 \frac { \Gamma(z/2) \zeta(z) } { 2 \pi^{z/2}  }
\end{equation}
where $\zeta(z) $ is the Riemann zeta function and $\Gamma(z)$  is the gamma function.   It follows that the zeros of the $\xi$ function  are identical to those of the $\zeta$ function.

The object of this paper is to investigate the properties of the $\xi$ function. The symmetric form of the $\xi$ function involves the integral of an infinite sum, and analysis is decidedly not easy.  Reformulation is possible by taking a finite sum of such terms, and then evaluating each integral; this makes sense if and only if the terms rapidly become negligible.  And if the integral can be written in an analytic form, then analysis becomes feasible. 

The plan of the paper is as follows.  In Section 2, the established properties of the $\xi$ function are described.  In Section 3, a reformulation of the $\xi$ function is carried out, and a validation exercise is undertaken. The objective of the paper occurs in Section 4, where the properties of the reformulated $\xi$ function are derived. A discussion and summary of these properties occurs in Section 5, while in Section 6, a level curve analysis is used to address the question as to whether these properties are informative in the context of the Riemann hypothesis.

\section{Established Properties}

 The $\xi$ function has the following established properties:

\begin{prop}It is an analytic function of complex $z$, and is thus holomorphic. 
\end{prop}
  Let $G$ be the region within the critical strip $0 \le \Re[z] \le 1$. Since $\xi$ is holomorphic in the region $G$ and is differentiable in the region $G$, then $\xi$ is infinitely differentiable and infinitely integrable in $G$.
\begin{prop} It satisfies a conjugation condition:
\[    
\xi(\overline{z}) = \overline{\xi(z)} 
\]
\end{prop}
\begin{prop} It is a symmetric function about the critical line $\Re[z] =0.5$. This symmetry is a consequence of the conjugation condition and the functional equation:
\[
   \xi(z) = \xi(1-z) 
\]
\end{prop}
A direct consequence of this symmetry is that $\Re[\xi]$ is an even function about the  critical line and $\Im [\xi]$ is an odd function about the  critical line.  Hence, on the critical line, $\xi$ is real.
 

\section{Reformulation}

The symmetric form of the $\xi$ function [3] is given by:
\begin{equation}
\label{eqn:symmetric}
\xi(z)	= -\frac{1}{z}+\frac{1}{z-1}+\int_1^\infty( x^{z/2-1}+x^{-(z+1)/2}) \; \omega(x)dx
\end{equation}
where 
\[\omega(x) = \sum_{n=1}^\infty e^{-\pi n^2 x} \]
  Note that Eqn. \ref{eqn:symmetric} is unchanged if $z$ is replaced by $1-z$, thus demonstrating Property 3. 

Consider the first term ($n=1$) in the expansion of $\omega(x)$, 
 and apply a change of variable: $t = \pi x$.
\begin{eqnarray*}
 \int_1^\infty x^{z/2-1} e^{-\pi x} dx  & = & \pi^{-z/2} \int_\pi^\infty t^{z/2-1} e^{-t} dt\\
                                         & = & \frac{\Gamma(z/2, \pi)}{\pi^{z/2}} 
\end{eqnarray*}
where $\Gamma(z/2, \pi)$ is the complex upper incomplete gamma function. Hence, define:
\begin{eqnarray}
v_0(z)   & = & -\frac{1}{z}+\frac{1}{z-1} \\
v_n(z)	& = & \frac{\Gamma(z/2, \pi n^2)}{(\pi n^2)^{z/2}} + \frac{ \Gamma((1-z)/2, \pi n^2)}{(\pi n^2)^{(1- z)/2}}  \hspace{.1in} \qquad n \in 1,2, \ldots \label{eqn:vn:gamma} \\
        & = &  e^{-\pi n^2} [U(1,1+z/2,\pi n^2) + U(1,1+(1-z)/2,\pi n^2)]
\end{eqnarray}
where $U(a,b,c)$ is the KummerU  confluent hypergeometric function [1].
Consequently, 
\begin{eqnarray}
  \xi_k(z) & = & \sum_{n=0}^k v_n(z) \label{eqn:xik} \\
  \xi(z)   & = & \sum_{n=0}^\infty v_n(z)   \label{eqn:xiinfty}
\end{eqnarray}

 Throughout, we use $z = x + iy$. In addition,  for historical consistency, the critical line is defined as $s = 0.5 + i \, t$. Hence, on the critical line:
\begin{eqnarray}
v_0(s)    & =  & -\Re[2/s]\\
         & = & \lim_{n \to 0}   \Re \left[ \frac{ \Gamma(  s/2 , \pi n^2)}{(\pi n^2)^{s/2}}\right] \\
         & = & \Re [ U(1,1+s/2,0)]
\end{eqnarray} 
and
\begin{eqnarray}
v_n(s)	& = & 2 \, \Re \left[ \frac{\Gamma(  s/2 , \pi n^2)}{(\pi n^2)^{s/2}}\right] \\
        & = & 2 \, \Re \, [e^{-\pi n^2} U(1,1+s/2,\pi n^2)]   \label{eqn:v_n} 
\end{eqnarray}

\begin{figure}[h] 
  \caption{Evaluation of $\xi_2(s)$ on the critical line }
  ~\newline
  \includegraphics[ width=3in, height = 2in]{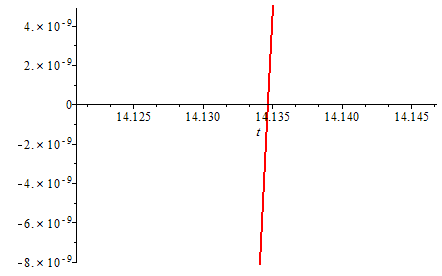}
  \label{fig:firstzero}
\end{figure}

To demonstrate this approach, Eqn. \ref{eqn:xik} is evaluated  in the region of the first non trivial zero using $k=2$. $\xi_2(s)$ is shown in Figure \ref{fig:firstzero}; the value of 14.1347 corresponds to the first non-trivial $\zeta$ zero.

\begin{table}[ht]
~\newline

~\newline

\caption{Values of $\xi_k(s)$ function on the critical line at a $\zeta$ zero}
\begin{tabular}{l l l l l l  }
\hline
 $t= $ & $14.134725$ & $21.022040$ & $25.010858$ & $30.424876$  & $32.935062$ \\
\hline
$\xi_1$  & -4.2240e-07  & -3.3657e-07  &  -2.8985e-07  & -2.3479e-07 & -2.1290e-07 \\
$\xi_2$  & -1.9938e-12  & -6.1334e-14  &  -2.9228e-14  & -2.8846e-14 & -2.7815e-14 \\
$\xi_3$  & -1.9594e-12  & -2.9001e-14  &   1.6977e-15  &  6.5564e-18 &  3.1076e-18  \\
$\xi$    & -1.9594e-12  & -2.9001e-14  &   1.6977e-15  &  6.5564e-18 &  3.1076e-18  \\
\hline 
\end{tabular}
\label{table:zero} 
\end{table}

~\newline

\begin{table}[ht]
~\newline

~\newline

\caption{Values of $\xi_k(s)$ function on the critical line}
\begin{tabular}{c l l l l l }
\hline
 $t= $ & $21.8$     & $26.2$     & $34.7$     & $41.6$      & $50.8$  \\
       & min        &  max       &  max        &  max       &  min     \\
\hline
$\xi_1$ &-3.5800e-08 & -2.7491e-07 & -1.9885e-07 & -1.5356e-07 &  -1.1180e-07  \\
$\xi_2$ &-3.0865e-08 &  1.9456e-09 &  3.4621e-12 & -1.9514e-14 &  -2.0909e-14  \\
$\xi_3$ &-3.0865e-08 &  1.9456e-09 &  3.4893e-12 &  4.8731e-15 &  6.1854e-18  \\
$\xi$   &-3.0865e-08 &  1.9457e-09 &  3.4893e-12  & 4.8731e-15 &  6.1854e-18  \\
\hline
$v_2$  & 3.2713e-07 & 2.7686e-07 &  1.9986e-07 &  1.5356e-07 & 1.1180e-07  \\
$v_3$  & 3.2069e-14 & 3.0483e-14 &  2.7146e-14 &  2.4387e-14 & 2.0903e-14  \\
$v_4$  & 5.5592e-24 & 5.4580e-24 &  5.2243e-24 &  5.0051e-24 &4.6859e-24  \\
\hline 
\end{tabular}
~\newline

~\newline

\label{table:xik} 
\end{table}

The accuracy of Eqn. \ref{eqn:xik} can be ascertained by evaluating $\xi_k(s)$ for $k \in \{1:3\}$ and comparing it to $\xi(s)$ at the first five $\zeta$ zeros{\footnote{Values evaluated using Maple, with Digits:=40.}.  These are shown in Table \ref{table:zero}.  $\xi_2$ is very close to $\xi$ for the first zero, while $\xi_3$ matches $\xi$ for the first five zeros.

Table \ref{table:xik} compares $\xi_k(s)$ to $\xi(s)$ on the critical line for a number of turning points.  For a small value of $t$, such as $21.8$, $\xi_2(s)$ accurately evaluates $\xi(s)$. From this table, it is clear that the effect of $v_3$ is insignificant, since it is about 6 orders of magnitude less than $\xi(s)$, which is why $\xi_2(s)$ provides all the accuracy that is needed.   However, as $t$ increases, $\xi(s)$ decreases exponentially,  so that by $t=34.7$ the effect of $v_3$ is no longer negligible.  At this $t$, it can be seen that $\xi_2(s)$ no longer provides a good enough estimate of $\xi(s)$, while $\xi_3(s)$ works well, and continues to do so for values of $t$ up to the mid 50s, for the same reason - the next term in the expansion, $v_4$, is of the order of 4e-24, which at $t=50.8$ is again 6 orders of magnitude smaller than $\xi(s)$  and so can be ignored.  The term ``epoch of order $k$" is used to signify the range of $t$ for which $\xi_k(s)$ accurately estimates $\xi(s)$.

\section{Properties derived from the Reformulated function}

Additional properties of $\xi_k(z)$, and hence $\xi(z)$, can be ascertained by evaluating  each term of $v_n(z)$ using Eqn. \ref{eqn:xik}.
 Because each term is relatively simple, this procedure is feasible.  In the following, given that $z = x+iy$,  we make the assumption that $y$ is large, and hence $x \ll y$.  
Define:
\begin{eqnarray}
       v^*_n(z) & = & \lim_{y \to \infty} v_n(z) \nonumber \\
      \phi_k(z) & = & \sum_{n=1}^k v_n^*(z)        \nonumber
\end{eqnarray}

\begin{prop}
Over the range $0 \le x \le 1$, $\Re [v_0^*(z)]$ is a negative convex function with a minimum at $x = 0.5$ for all $y$.
\end{prop}
\textbf{Proof:}
By definition: 
\begin{eqnarray}     
           v_0(z)       & = & -1/z + 1/(z-1) \nonumber \\
          \Re[v_0(z)]   & = & -x/(x^2+y^2) + (x-1)/((x-1)^2+y^2) \nonumber \\
          \Re[v_0^*(z)] & \approx & (x (x-1) -y^2)/y^4 \nonumber \\
                        & = & -(y^2 - r^2)/y^4  \label{eqn:v0}
\end{eqnarray}
where $ r=x-0.5$.  Thus for any given $y$, $\Re [v^*_0(z)]$ in the $x$ plane is a negative convex quadratic function about $x = 0.5$, while in the $y$ plane, $\Re [v^*_0(z)] $ asymptotically approaches zero from below.  
\qedsymbol

\begin{prop}
 On the critical line, $\Re [v_0^*(s)] $ is negative throughout, has a minimum value of -4 at $t=0$, and increases as a concave function asymptotically reaching zero as $t \to \infty.$ 
\end{prop}
\textbf{Proof:}
This follows directly from the definition of $v_0^*$ and the assymptotic behaviour depicted in Eqn. \ref{eqn:v0}.
\qedsymbol

\begin{prop}
Over the range $0 \le x \le 1$, $\Re [v_n^*(z)]$ is a positive concave function with a maximum at $x = 0.5$ for all $n$ and for all $y$. 
\end{prop}

\textbf{Proof:}
Consider first $v_1(z)$.  Recall from Eqn. \ref{eqn:vn:gamma}  that:
\begin{equation}
v_1(z)	=  \frac{\Gamma(z/2, \pi)}{(\pi)^{z/2}} + \frac{ \Gamma((1-z)/2, \pi)}{(\pi)^{(1- z)/2}}
\end{equation}
The lower incomplete gamma function has a power series expansion:
\begin{equation}
          \gamma(z, x) = x^z \Gamma(z) e^{-x} \sum_{k=0}^\infty \frac {x^k}{\Gamma(z+k+1)}
\end{equation}
Thus using the fact that $\Gamma(z+1) = z \Gamma(z)$, $\Gamma(z,\pi) = \Gamma(z) -\gamma(z,\pi)$,  and using just the first two terms of the expansion, since the terms converge extremely rapidly, we have:
\begin{eqnarray*}
   \frac{\Gamma(z/2, \pi)}{(\pi)^{z/2}} & = & \frac {\Gamma(z/2)}{\pi^{z/2}} -2 e^{-\pi}\left[ \frac {1}{z} + \frac {2 \pi}{z(z+2)} \right] \\
   \frac{\Gamma((1-z)/2, \pi)}{(\pi)^{(1-z)/2}} & = & \frac {\Gamma(-z/2+.5)}{\pi^{(1-z)/2}} -2 e^{-\pi}\left[ \frac {1}{1-z} + \frac {2 \pi}{(1-z)(3-z)} \right] 
\end{eqnarray*}
For large values of $y$, the term $\Gamma(z)/\pi^z$ rapidly becomes insignificant - for $y=100$, this term has a real component of the order of 1e-68. This is also true for larger $n$; for example, with $y=1000$ and $n = 4$, $\Gamma(z/2)/(n^2 \pi)^{z/2}$ has a real component of the order of 1e-342.  Thus:
\begin{eqnarray}
       v_1(z) & \approx &  -2 e^{-\pi}\left[ \frac {1}{z} + \frac {2 \pi}{z(z+2)} + \frac {1}{1-z} + \frac {2 \pi}{(1-z)(3-z)} \right] \label{eqn:v1:exp} \\
             & = & 2 e^{-\pi} \frac {(4 \pi-1)(z-z^2)+6(1+\pi)}{z (z+2) (z-1) (z-3)} \nonumber \\ 
         \Re[v_1^*(z)]    & = & 2 e^{-\pi} \frac {(4 \pi-1)(x-x^2+y^2)}{y^4} \nonumber \\
               & = & 2 e^{-\pi} (4 \pi-1)(y^2-r^2)/y^4
\end{eqnarray}
where $ r=x-0.5$.  
For the general case, we have:
\begin{eqnarray}  
      v_n(z)           & \approx &  -2 e^{-n^2 \pi}\left[ \frac {1}{z} + \frac {2 n^2 \pi}{z(z+2)} + \frac {1}{1-z} + \frac {2 n^2 \pi}{(1-z)(3-z)} \right] \nonumber \\
         \Re[v_n^*(z)]     & = & 2 e^{-n^2\pi} (4 n^2 \pi-1)(y^2-r^2)/y^4  \label{eqn:vn*}
\end{eqnarray}
 Thus for a given $y$, $\Re [v^*_n(z)]$ in the $x$ plane is a positive concave quadratic function about $x = 0.5$,  while in the $y$ plane, $\Re [v^*_n(z)]$ asymptotically approaches zero from above. 
\qedsymbol

Since the sum of concave functions is concave, $ \Re [\phi_k(z)]$ in the $x$ plane is a positive concave function with a maximum at $x = 0.5$ for all $y$.

A comparison of Eqns.  \ref{eqn:v_n} and  \ref{eqn:vn*}  shows that: 
\begin{equation}
\lim_{t \to \infty}  \Re[U(1,1+s/2,\pi n^2)] =   (4 n^2 \pi-1)/t^2 
\label{eqn:kummeru}
\end{equation}

\begin{prop}
On the critical line, $\Re [v_n^*(s) \mid n>0]$ is positive throughout, declines exponentially,  and asymptotically tends to zero as $t \to \infty$.
\end{prop}

\textbf{Proof:}
This follows directly from  Eqn. \ref{eqn:vn*}.
\qedsymbol

\begin{prop}
   On the critical line, for any $k$, 
\[    \lim_{t \to \infty} \xi_k(s) \to -v^*_{k+1}(s)  \]
\end{prop}

\textbf{Proof:}
On the critical line, using Eqns.  \ref{eqn:v0} and  \ref{eqn:vn*},  $v_n^*(s) = \lambda_n/t^2$ for large $t$ and where
\begin{eqnarray}
   \lambda_0 & =  &  - 1    \\
   \lambda_n & =  & 2 e^{-n^2\pi} (4 n^2  \pi - 1) \hspace{.5in} \qquad n \in 1,2, \ldots  \label{eqn:lambdan}
\end{eqnarray}
$\lambda_n$ is a remarkable series because of the following two properties:
\begin{eqnarray}
 \sum_{n=0}^k \lambda_n  &  \approx & -\lambda_{k+1} \label{lem:lambda} \\   
 \sum_{n=0}^\infty \lambda_n & =& 0 
\end{eqnarray}
Hence:
\begin{eqnarray*}
  \lim_{t \to \infty} \xi_k(s) & = & \sum_{n=0}^k v^*_n(s) \\
                               & = & \sum_{n=0}^k \lambda_n(s)/t^2 \\
                               & = & -\lambda_{k+1}/t^2 \\
                               & = & -v^*_{k+1}(s) \hspace{.1in} 
\end{eqnarray*}
 \qedsymbol

 The properties of $\lambda_n$ are illustrated in  Table \ref{table:lamda}.  The first column shows $\lambda_n$; this follows directly from the definition in Eqn. \ref{eqn:lambdan}, and is responsible for the damping of the $\xi$ function.
 The second column shows the sum of the first $k$ terms of $\lambda_n$, and clearly illustrates
 Eqn. \ref{lem:lambda}. 
The limiting behaviour of $\lambda_k $ and $\xi_k$ respectively can be
expressed as  ratios; these are given by:
\begin{eqnarray*}
  p_k    & = &  \frac {\sum_{n=0}^k \lambda_n} { -\lambda_{k+1}} -1 \\
  q_k(s) & = & \lim_{t \to \infty} \frac {\sum_{n=0}^k v _n} { -v_{k+1}(s)} -1        \\      
         & = &  \lim_{t \to \infty} \frac {\xi_k} { -v_{k+1}(s)} -1      
\end{eqnarray*}
The third and fourth columns display these ratios where, for the fourth column, $s^* = .5+1000000i$.  The fact that these two columns are identical provides a numerical confirmation of this property.

\begin{table}[ht]
\caption{Values of $\lambda_n$ }
\begin{tabular}{c l l l l  }
\hline
$n..k$  & $~~~~~\lambda_n $ & ~~~~~$\sum_{n=0}^k \lambda_n$ & $ p_k $ & $q_k(s^*) $   \\
\hline
0  & -1.0000e-00  & -1.0000e-00 &  .3436e-03  & .3437e-03  \\
1  &   .9996e-00  & -.3436e-03  &  .3429e-06  & .3429e-06 \\
2  &   .3436e-03  & -.1178e-09  &  .5023e-09  & .5023e-09 \\
3  &   .1178e-09  & -.5918e-19  &  .8226e-12  & .8226e-12 \\
4  &   .5918e-19  & -.4868e-31  &  .1415e-14  & .1415e-14 \\
5  &   .4868e-31  & -.6887e-46  &  .2496e-17  & .2496e-17 \\
\hline 
\end{tabular}
\label{table:lamda} 
\end{table}

\begin{prop}
Over the range $0 \le x \le 1, \Im [v_0^*(z)]$ is a negatively sloped straight line with a value of zero at $x = 0.5$ for all positive $y$. 
\end{prop}

\textbf{Proof:}
By definition:  
\begin{eqnarray}     
      v_0(z)         & = &  -1/z + 1/(z-1) \nonumber \\
      \Im[ v_0(z)]   & = & y/(x^2+y^2) - y/((x-1)^2+y^2) \nonumber \\
      \Im[ v_0^*(z)] & \approx & ( (2x-1))/y^3 \nonumber \\
                     & = & - 2 r/y^3
\end{eqnarray}
where $ r=x-0.5$.  Thus for a fixed positive $y$, $\Im [v^*_0(z)]$ in the $x$ plane is a negatively sloped straight line with a value of zero at $x = 0.5$, while in the $y$ plane, $\Im [v^*_0(z)]$ asymptotically approaches zero from above or below, depending on the sign of $r$.
\qedsymbol

\begin{prop}
Over the range $0 \le x \le 1, \Im (v_n^*(z))$ is a positively sloped straight line with a value of zero at $x = 0.5$ for all $n$ and for all positive $y$. 
\end{prop}

\textbf{Proof:}
Recall Eqn: \ref{eqn:v1:exp}:
\begin{eqnarray}
   v_1(z) & \approx &  -2 e^{-\pi}\left[ \frac {1}{z} + \frac {2 \pi}{z(z+2)} + \frac {1}{1-z} + \frac {2 \pi}{(1-z)(3-z)} \right] \nonumber \\
                 & = & 2 e^{-\pi} \frac {(4 \pi-1)(z-z^2)+6(1+\pi)}{z (z+2) (z-1) (z-3)} \nonumber \\ 
   \Im[ v_1^*(z)] & = & 2 e^{-\pi} \frac {(4 \pi-1) (2x-1)}{y^3} \nonumber \\
                 & = & 4 r e^{-\pi} (4 \pi-1)/y^3
\end{eqnarray}
where $ r=x-0.5$.  
For the general case, we have:
\begin{eqnarray}
    v_n(z) & \approx &  -2 e^{-n^2 \pi}\left[ \frac {1}{z} + \frac {2 n^2 \pi}{z(z+2)} + \frac {1}{1-z} + \frac {2 n^2 \pi}{(1-z)(3-z)} \right] \nonumber \\
        \Im[ v_n^*(z)]     & = & 4 r e^{-n^2\pi} (4 n^2 \pi-1)/y^3
\end{eqnarray}
Thus for any given $y$, $\Im [v^*_n(z)]$ in the $x$ plane is a positively sloped straight line with a value of zero at $x = 0.5$, while in the $y$ plane, $\Im [v^*_n(z)]$  asymptotically approaches zero from above or below, depending on the sign of $r$.
\qedsymbol

\begin{prop}
The asymptotic periodicity $T(z/2,\pi)$ of the upper incomplete gamma function $\Gamma(z/2,\pi) $ is:
\[ \lim_{y \to \infty} T(z/2,\pi) = \frac{4 \pi}{\ln(\pi)} \]
\end{prop}

\textbf{Proof:}
From the proof of Property 6, the power series expansion of the upper  incomplete gamma function is:
 \[  \Gamma(z/2, \pi)  = -2 \pi^{z/2} e^{-\pi}\left[ \frac {1}{z} + \frac {2 \pi}{z(z+2)} + \cdots \right] \]
The only periodic component is in the term $\pi^{z/2}$, since the remaining terms only affect the amplitude. Thus it follows that the asymptotic periodicity  of $\Gamma(z/2,\pi) $ is the same as the periodicity of  $\pi^{z/2}$. Since:
\[ \pi^{z/2} = \pi^{x/2}[\cos( \;  \ln(\pi) y/2 \; )+ i \sin( \;  \ln(\pi) y/2 \; )] \] 
the periodicity is $4 \pi/\ln(\pi)$. \qedsymbol

 This proof can be made stronger, as shown below:

\begin{prop}
The asymptotic periodicity of the the upper incomplete gamma function outside the critical strip is the same as the periodicity on the critical line.
\end{prop}

\textbf{Proof:}
Given the standard Kummer relationship [1] 
\begin{equation}
       c U(a,b+1,c) = (b-a) U(a,b,c) + U(a-1,b,c)  \nonumber
\end{equation}
Let $a = 1$, $ b=1+s$, and $c= \pi $, and given $U(0, b, c) = 1 $,
\begin{eqnarray}
       \pi U(1,b+1,\pi) & = &  (b-1) U(1,b, \pi) + 1 \nonumber \\
       \pi U(1,s+2,\pi) & = &   s \, U(1,s+1,\pi) + 1 \nonumber\\
       \pi e^{\pi} \pi^{-(s+1)} \Gamma(s+1,\pi) & = & s \, e^{\pi} \pi^{-s} \Gamma(s,\pi) + 1  \nonumber \\
       \Gamma(s+1,\pi) & = &  s \,\Gamma(s,\pi) +  e^{-\pi} \pi^{s} \label{eqn:s1}
\end{eqnarray}
Let $T(s,\pi)= 2 \pi /\ln(\pi) $ be the  asymptotic periodicity of the incomplete gamma function $\Gamma(s,\pi)$ on the critical line.
Then from Eqn \ref{eqn:s1}:
\begin{equation}
     \Gamma(s+iT+1,\pi)  =   (s+iT) \,\Gamma(s+iT,\pi) +  e^{-\pi} \pi^{s+iT} \label{eqn:s1t1}
\end{equation}
Combining Eqns \ref{eqn:s1} and \ref{eqn:s1t1}:
\begin{eqnarray*}
 \Gamma(s+iT+1,\pi)- \Gamma(s+1,\pi) &=& (s+iT) \, \Gamma(s+iT,\pi) - s \,\Gamma(s,\pi) \\
                                      &  &  + \; e^{-\pi} \pi^{s} (\pi^{iT}-1) 
\end{eqnarray*}
Since $\lim_{y \to \infty} (s+iT) \, \Gamma(s+iT,\pi) - s \,\Gamma(s,\pi) = 0$ and $(\pi^{iT}-1) = 0$, we have $\lim_{y \to \infty} \Gamma(s+iT+1,\pi)- \Gamma(s+1,\pi) = 0$, thus proving that the asymptotic periodicity of the upper incomplete gamma function at $s+1$ outside the critical strip is the same as the periodicity on the critical line. \qedsymbol


\begin{prop}
On the critical line, the maximum (minimum) of the $\xi$ function in the y plane occurs at the intersection of the real level curve $\Im [\xi(z)] = 0$.
\end{prop}

\textbf{Proof:}
On the critical line, $\xi(z)$ is real.  Thus on the critical line, 
  $ \frac{\partial  \Re [\xi(z)]}{\partial y} = 0$
corresponds  to the top of a ridge or the bottom of a valley of $\Re [\xi(z)]$.  
Similarly, on the critical line,  $\frac{\partial \Re [\xi(z)]}{\partial x} =0$
 because 
$\Re [\xi(z)]$ is even about the critical line.  Since the Laplace condition is satisfied, this is a  saddle point.  From the Cauchy-Riemann equations,
 $ \frac{\partial  \, \Im [\xi(z)]}{\partial y} = 0$
 and 
 $\frac{\partial  \, \Im [\xi(z])}{\partial x} = 0$. 
The former follows immediately from the fact that on the critical line, $\Im [\xi(s)] =0$, while the latter situation occurs when $\Im [\xi(z)]$ is constant and has a value of zero;  this  occurs when
 $\Im [\xi(z)] =0$  - that is the real level curve.  \qedsymbol

\section{Discussion}
 
\subsection{Damping}
The $\xi(s)$ function decreases exponentially with $t$. Thus in the range $t=21.8$ to $50.8$, $\xi(s)$ decreases by 10 orders of magnitude.  This is a consequence of two components.  The first is the dampening effect of  $v_n(s) $ because of the term $1/t^2$ shown in Eqn. \ref{eqn:kummeru}. The second, and much stronger, dampening effect is a consequence of Property 8, whereby the  function $\xi_k(s)$ tends to the  limit $-v_{k+1}(s)$.  Both these effects can be seen in Table 2.  Thus $\xi(s)$ has a very short decay time.

\subsection{Epoch}
 The term ``epoch of order $k$" is used to signify the range of $t$ for which $\xi_k(s)$ accurately estimates $\xi(s)$.  When  $t$ increases sufficiently such that $\xi_k(s)$ is sufficiently small and is affected (at the level of accuracy required) by $v_{k+1}(s)$, then the next epoch is required to accurately estimate $\xi(s)$.  In ordinary terms, we would state that the number of terms ($k$) required should be increased.  In this context, the implication is that the epochs are similar, in that the same behaviour is exhibited for each epoch, and thus the properties associated with one epoch is true for all epochs.

\subsection{Periodicity}
The $\xi$ function is not periodic, since there is no value $T$ such that $\xi(z+T)= \xi(z)$. 
Rather it is almost periodic, defined loosely in that a value is repeated at almost equally spaced intervals.  To show this, consider $\xi_k(z)$ which is composed of the sum of $v_n(z)$.  Recall from Eqn. 14 that   
$v_1(z)$ consists of two terms; each term being the ratio of an incomplete gamma function in the numerator, and $\pi$ raised to a complex power in the denominator.  The periodicity of the denominator ($\pi^{z/2}$) is $4 \pi /\ln(\pi) = 10.977586$,  while from Property 11, the periodicity of the numerator approaches the same value asymptotically.  Thus these two terms are incommensurable.  Hence, since the product of incommensurable periodic functions results in a function that is almost periodic, it follows that $v_1(z) $ is almost periodic, and similarly for $v_k(z)$.  Since the sum of a finite number of almost periodic functions is itself an almost periodic function, we have $\xi_k(z)$ is an almost periodic function. 

As $y$ increases, so the periodicity of the numerator asymptotically approaches the periodicity of the denominator.  Thus since the periodicity of $\pi^z$ is independent of $x$, so too the asymptotic periodicity of $\Gamma(z,\pi)$ is also independent of $x$.  This is a strong result, since the almost periodic nature of the $\xi_k$ function is a consequence of the almost periodic nature of the $v_n$ components.  In particular, focusing on $\Re [\xi(z)]$, we have the following assertion:   
\begin{assert}
For every real ridge (valley) of the $\xi$ function that exists outside the critical strip, there will be a corresponding real ridge (valley) that crosses the critical line. 
\end{assert}

\comm{

\footnote{Even though the two terms are incommensurable, they have almost identical periods, so $v_1(z)$ would be expected to have an approximate period of about the same value. An order of magnitude estimate of the  initial period of $\xi_1(z)$ can be evaluated as the difference between the third zero and the first zero - $ 25.0109-14.1347 = 10.8761$.}
While in no way a proof, it is instructive to look at the average spacing for the first epoch.  The first epoch is easy, since there is only one periodic component ($v_1)$, so there are no concerns over superposition from the sum of periodic components.   An order of magnitude estimate of the  initial period of $\xi$ can be evaluated as the difference between the third zero and the first zero - $ 25.0109-14.1347 = 10.8761$.  The mid value of that range is 19.5728,

\begin{table}[ht]
\caption{Comparison period estimation)}
\begin{tabular}{ l l l }

\hline
Actual zero range                 & 25.0109-14.1347       &10.8761   \\
Asymptotic period of $v_1$        &  $ 4 \pi /\ln(\pi)$   & 10.9776   \\
Twice average space between zeros &  $ 4 \pi / \ln(\bar{t}/2 \pi)$ & 11.058 \\ 
\hline 
\end{tabular}
\label{table:period} 
\end{table}
}  

\subsection{Oscillation}
On the critical line, Eqn. 12 depicts an under-damped oscillating function, $v_n(s)$.  $U(1,b,z)$ is related to $M(1,b,z)$
by [1]:
\begin{equation}
  U(1,b,z) = \Gamma(b-1) [z^{1-b} e^z - M(1,b,z)/\Gamma(b)] 
\end{equation}
where $M(1,b,z)$ is the KummerM confluent hypergeometric function.

In general,  $M(a,b,z)$ is entire in $z$ and $a$, and is a meromorphic function in $b$.  Consequently, $U(1,1+s/2, \pi n^2)$ shares the same properties, and can be represented by a complex trigonometric polynomial on $\Im[1+s/2] = t/2$. Thus, $v_n(s) $  oscillates, and consequently $\xi_k(s)$ also oscillates, and these oscillations are damped since $\xi_k(s)$ has a very short decay time.

\begin{figure}[h] 
 \caption{$\xi_1(s)$ and $-v_2$ on the critical line }
  ~\newline
  \includegraphics[ width=3in,  keepaspectratio]{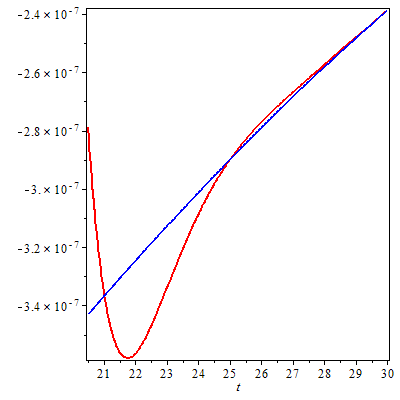}
  \label{fig:xi1}
\end{figure}

\begin{figure}[h] 
 \centering
  \caption{$\xi_2(s)$ and $-v_3$ on the critical line }
  ~\newline
 \includegraphics[ width=3in, keepaspectratio]{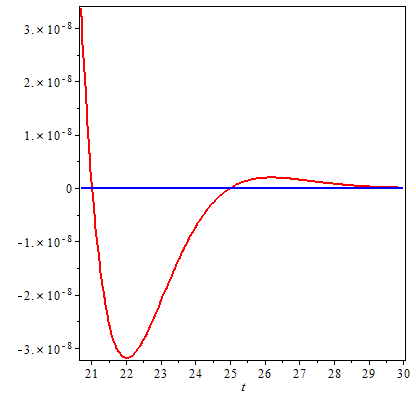}

  \label{fig:xi2}
\end{figure}

   From the discussion on periodicity, the  $\xi$ function is an almost periodic function, and the behaviour of the $\xi$ function can be determined by an examination of a particular epoch.  Each epoch of value $k$ is characterized by an almost periodic oscillating function, and is terminated when $v_{k+1}(s)$ becomes significant.  

Consider the first epoch, part of which is depicted in Figure 2. $\xi_1(s)$ is shown in red and $-v_2(s)$ in blue.\footnote{By displaying $v_2$ as negative, we can focus on the interaction between $\xi_1$ and $v_2$.} From Property 8, we have 
\begin{equation}
\lim \limits_{t \to \infty} \xi_1(s) + v_{2}(s) = 0 
\label{eqn:limit}
\end{equation}
$\xi_1(s)$ is  both a damped and an almost periodic oscillating function. It has a non linear asymptote $-v_2(s)$, and thus has maxima greater than $-v_2(s)$ and minima less than $-v_2(s)$. Because of its short decay time, oscillations die out very quickly, as can be seen from Figure 2.

Clearly, over the range shown, and as can be seen in Table \ref{table:xik}, $\xi_1(s)$ is not a good measure of $\xi(s)$ because $v_2(s)$ is of the same order of magnitude. Consequently, consider the  next epoch, $\xi_2(s)$ and $-v_3(s)$; these are depicted in Figure 3.
 
Over the same range, $v_3(s)$ is a number of orders of magnitude less than $\xi_2(s)$, and can be treated as zero.  $\xi_2(s)$ and $\xi(s)$ are effectively identical, and thus over this epoch and from Property 8 we have :
\begin{equation}
\lim \limits_{t \to \infty} \xi_2(s)  = -v_3(s) \approx 0 
\end{equation} 
The asymptote for $\xi_2(s)$ for this epoch is thus zero, and thus has maxima greater than $-v_3(s)$ and minima less than $-v_3(s)$; this implies the following assertion:
\begin{assert}
On the critical line, the maxima of the $\xi$ function are positive, and the minima are negative. 
\end{assert}

\subsection{Summary}Thus the behaviour of the $\xi$ function can be explained by the behavior in each epoch. This  consists of a damped almost periodic oscillation around a limiting nonlinear asymptote  that is a number of orders of magnitude less than the turning points within the epoch.

\section{Level Curve Analysis}

This section considers  the question as to whether the properties of the $\xi$ function described above can be used to address the Reimann hypothesis. 

The real part of the $\xi$ function can be described as a sequence of ridges and valleys - the ridges descending and the valleys rising to the critical line. Formally, assuming RH, a ridge results in a maximum on the critical line, and a valley results in a minimum. 
 
A non trivial  zero of the $\xi$  function occurs at the intersection of the two level curves The term [4] ``real level curve" describes the line  $\Im [\xi(z)]= 0$ and ``imaginary level curve" the line   $\Re [\xi(z)]= 0$.   There are only four possible scenarios, which are depicted in Figure \ref{fig:Zeros}:

\begin{figure}[h] 
  \includegraphics[ width=6in, height = 4in, keepaspectratio]{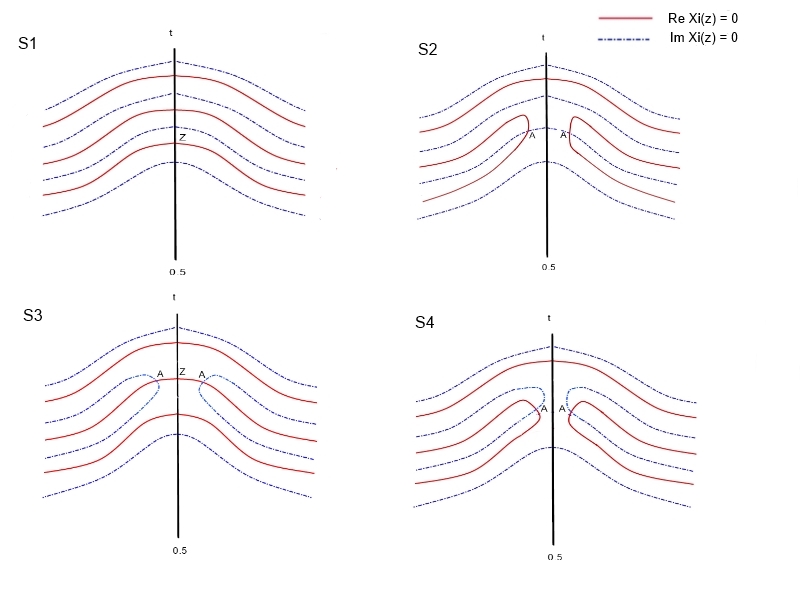}
  \caption{Scenarios for the zeros of the $ \xi$ function }
  \label{fig:Zeros}
\end{figure}

The easiest way of generating zeros off the critical line is to envisage the ordered set of ridges\footnote{The discussion is in terms of the ridge. To simplify the discussion, neighboring minima are assumed to exhibit strictly negative values. The equivalent analysis can obviously be applied to the valley.} in decreasing magnitude, followed by the set of ridges that die out before crossing the critical line. This set is categorized as follows:
\begin{itemize} 
  \item $R_1$. The ridges that cross the critical line and have a critical point that is positive, and hence satisfy the RH. 
  \item $R_2$. The ridges that cross the critical line and have a critical point that is negative, and hence do not satisfy the RH. 
  \item $R_3$. The ridges that die out before crossing the critical line, and hence also do not satisfy the RH.
\end{itemize}  
Thus an $R_2$ ridge differs from an $R_1$ ridge only in the sign of the critical point.

~\newline

\textbf{Scenario: S1}
\newline
 This scenario, which occurs under R1, is characterized by the intersection of the level curve  $\Re [\xi(z)]= 0$ and the critical line - this creates a single zero (Z) on the critical line. The Riemann hypothesis is that all zeros fall on the critical  line [2], and consequently if the hypothesis is true, then only S1 is viable. 

\begin{figure}[h] 
  \centering
\includegraphics[width=4in, height = 3.5in, keepaspectratio]{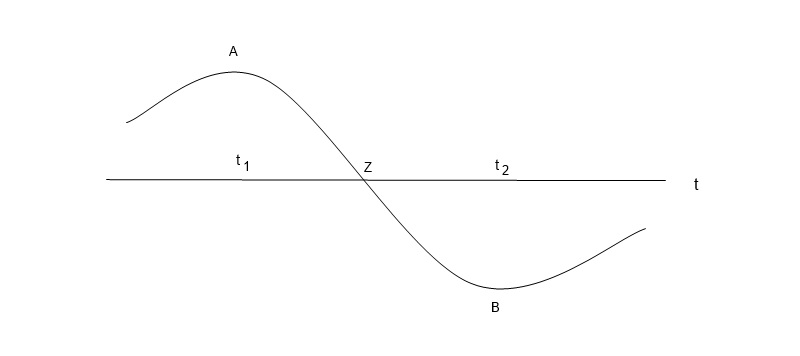}
  \caption{S1 - Elevation view of $\Re [\xi(s)]$  }
  \label{fig:s1}
\end{figure}

Under R1 a ridge crosses the critical line with a positive value; this is shown in Figure \ref{fig:s1}.  The panel depicts $ \xi(s)$ on the critical line. $\xi(s)$ is real, and has critical points  A (a maximum)  at $t_1$, and B (a minimum) at $t_2$ in the y plane.  Thus from Property 13,
the real level curve, $\Im [\xi(z)] = 0$, must intersect the critical line at $t_1$ and $t_2$.  Since A is positive and B is negative, there must be a point Z between them on the critical line where $\Re [\xi(s]) = 0 $, which will lie on the imaginary level curve $\Re [\xi(z)] = 0$.  Thus Z is a zero of the $\xi$ function, and thus satisfies the RH. 

\textbf{Scenario: S2}
\newline
This scenario occurs under R2. Under R2, a ridge crosses the critical line with a negative value; this is shown in Figure \ref{fig:s2}.
Under Property 13,  the real level curve, $\Im [\xi(z)] = 0$, must still intersect the critical line at $t_1$ and $t_2$, irrespective as to the signs of A and B. As in S1, there is a maximum at A and a minimum at B, but unlike S1, the maximum and minimum have the same sign, and thus there is no point between them that has a zero.  Consequently, the zeros would have to occur off the critical line, as depicted in S2 of Figure 4.  However, this situation violates Assertion 2, whereby maxima are positive and minima are negative.  Hence S2 is not viable.

\begin{figure}[h] 
  \centering
  \includegraphics[width=4in, height = 3.5in, keepaspectratio]{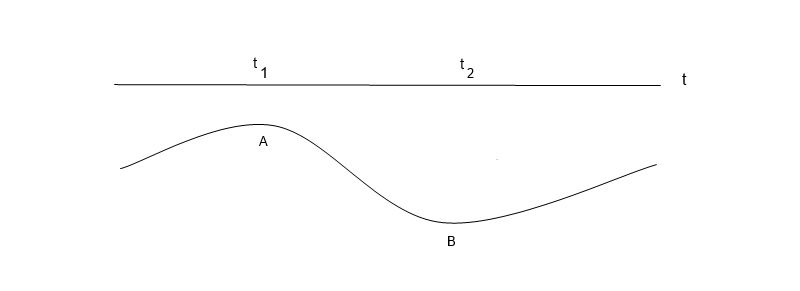}
  \caption{S2 - Elevation view of $\Re [\xi(s)]$  }
  \label{fig:s2}
\end{figure}

\textbf{Scenario: S3}
\newline
Next, consider scenario S3. As can be seen in Figure \ref{fig:Zeros}, there are a pair of zeros shown at (A), as well as a zero (Z) on the critical line.  All three zeros lie on the level curve $\Re [\xi(z)] = 0$.  Since the harmonic conjugate of the level curve of a meromorphic function is monotonic,  $\Re [\xi(z)] = 0$ is required to be monotonic on  each side of the critical line.  Clearly this is not the case, and hence S3 is not viable.

\textbf{Scenario: S4}
\newline
Finally, in scenario S4, the ridge is of type $R_3$ - it dies out before reaching the critical line.  The situation is shown in  Figure \ref{fig:s4}.   Using an elevation view along the x plane, a ridge descends as it approaches the critical line.  However, unlike $R_1$ or $R_2$, the ridge doesn't reach the critical line, but dies out before hand in the surrounding valley which ascends to a maximum on the critical line (C).  Since the downward moving ridge joins the upward moving valley, there will be a minimum, shown at point A, and a point of inflection at point B.
 However, this situation violates Assertion 1, whereby every ridge (valley) outside the critical strip must cross the critical line.  Hence S4 is also not viable.

\begin{figure}[h] 
  \centering
  \includegraphics[width=4in, height = 3.5in, keepaspectratio]{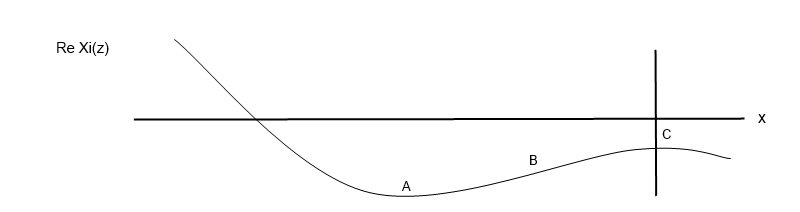}
  \caption{S4 - Elevation view of $\Re [\xi(s)]$ on the real axis. }
  \label{fig:s4}
\end{figure}

\section{Conclusion}
This paper has demonstrated that the $\xi$ function can be represented as the weighted sum of incomplete gamma functions, and this formulation has permitted a number of properties of the $\xi$ function to be exposed. It has also been demonstrated that the existence of a zero of the $\xi$ function, (and hence by extension of the $\zeta$ function) off the critical line is not compatible with these properties, thus validating the Riemann hypothesis.


\section*{Acknowledgment}
I would like to thank Marco Bertola, Ken Davidson, Lubo\u{s} Motl and Terence 
Tao for their  feedback and patience.

This work is dedicated to the memory of my grandfather, Lou Jaffe.


\begin{thebibliography}{10}


\bibitem {1} M. Abramowitz and I.A. Stegun, (Eds.), (1972). \emph{Handbook of
Mathematical Functions with Formulas, Graphs, and Mathematical Tables}, New
York: Dover.

\bibitem {2} J.B. Conrey, (2003). ``The Riemann Hypothesis'' \emph{Notice of the AMS}, 50, 341-353.

\bibitem {3} H.M. Edwards, (1974). \emph{Riemann's Zeta Function}, Academic Press, New York.

\bibitem{4} J.K. Hill and R.M. Wilson, (2004).``X-rays of the Riemann Zeta and Xi functions'', \emph{www.numbertheory.org/pdfs/xray.pdf}.


\end{thebibliography}

~\newline 
\bibliographystyle{amsplain}

\end{document}